\documentclass[11pt]{article}
\usepackage{amssymb,latexsym}
\usepackage[dvips]{graphicx} 
\usepackage{url}
\usepackage{picture}
\usepackage{color}

\newif\iftemp \tempfalse 
\temptrue

 \catcode`@=11\def\@oddhead{\hbox{}\hfil\rm\thepage}\def\@oddfoot{}
 \def\@evenhead{\hbox{}\hfil\rm\thepage}\def\@evenfoot{}
 \@twosidetrue \catcode`@=12

\iftemp
\DeclareFixedFont{\itshape}{OT1}{cmr}{m}{it}{11}

\fi
\newtheorem{prp}{Proposition}
\newtheorem{lem}[prp]{Lemma}\newtheorem{thm}[prp]{Theorem}

\newtheorem{cor}[prp]{Corollary}
\newenvironment{prf}{\begin{trivlist}\item[\emph{Proof.}]}{\end{trivlist}
  \medskip\par}
\newenvironment{rem}{\begin{trivlist}\item[\emph{Remarks.}]}{\end{trivlist}
  \medskip\par}
\def\prpb{\begin{prp}}\def\prpe{\end{prp}}
\def\lemb{\begin{lem}}\def\leme{\end{lem}}
\def\thmb{\begin{thm}}\def\thme{\end{thm}}
\def\corb{\begin{cor}}\def\core{\end{cor}}
\def\prfb{\begin{prf}}\def\prfe{\end{prf}}
\def\remb{\begin{rem}}\def\reme{\end{rem}}
\def\prpa#1{\label{p:#1}}\def\prpu#1{Proposition~\ref{p:#1}}

\def\thma#1{\label{t:#1}}\def\thmu#1{Theorem~\ref{t:#1}}

\def\seca#1{\label{s:#1}}\def\secu#1{Section~\ref{s:#1}}

\def\itmb{\begin{enumerate}}\def\itme{\end{enumerate}}
\def\itdb{\begin{itemize}}\def\itde{\end{itemize}}
\def\ittb{\begin{description}}\def\itte{\end{description}}
\def\eqnb{\begin{equation}}\def\eqne{\end{equation}}
\def\arrb#1{\begin{array}{#1}}\def\arre{\end{array}}
\def\tabb#1{\par\noindent\begin{tabular}{#1}}
\def\tabe{\end{tabular}\par\noindent}
\def\eqna#1{\label{e:#1}}\def\eqnu#1{(\ref{e:#1})}

\def\QED{\relax\ifmmode\let\@tempa\relax\ifcase\@eqcnt\def\@tempa{& & &}\or
  \def\@tempa{& &}\else\def\@tempa{&}\fi\@tempa $\Box$ \else\hfill $\Box$ \fi}
\def\DDD{\relax\ifmmode\let\@tempa\relax\ifcase\@eqcnt\def\@tempa{& & &}\or
 \def\@tempa{& &}\else\def\@tempa{&}\fi\@tempa $\Diamond$
 \else\hfill $\Diamond$ \fi}

\def\Rom#1{\uppercase\expandafter{\romannumeral#1}}
\def\dsp{\displaystyle}

\def\Ccomb#1#2{\setbox0=\hbox{$\displaystyle\mathrm{C}$}\setbox1=\hbox{%
$\scriptstyle #1$}\kern \wd1{\mathrm{C}}_{\kern -1.05\wd0\kern -0.99\wd1{#1}
 \kern 1.15\wd0{#2}}}

\def\clvec#1#2#3{\def\clvecone{#3}\left(\arrb{c} \dsp #1\\ \dsp #2
 \ifx\clvecone\empty\else\\ \dsp #3\fi\arre\right)}

\def\le{\leqq} \def\leq{\leqq} \def\geq{\geqq}

\def\twreals{{\mathbb R^2}}

\def\integers{{\mathbb Z}}\def\pintegers{{\mathbb Z}_+}
\def\nintegers{{\mathbb N}}

\def\prb#1{\def\prbone{#1}
  \ifx\prbone\empty{\mathrm{P}}\else{\mathrm{P[\;}}#1{\mathrm{\;]}}\fi}
\def\prbseq#1#2{\def\prbseqone{#2}
  \ifx\prbseqone\empty{\mathrm{P}}_{#1}\ignorespaces
  \else{\mathrm{P}}_{#1}{\mathrm{[\;}}#2{\mathrm{\;]}}\fi}

\def\EEseq#1#2{\def\EEseqone{#2}
  \ifx\EEseqone\empty{\mathrm{E}}_{#1}\else
 {\mathrm{E}}_{#1}{\dsp\mathrm{[\;}}#2{\mathrm{\;]}}\fi}

\def\VVseq#1#2{\def\VVseqone{#2}
  \ifx\VVseqone\empty{\matrm{V}}_{#1}\else
 {\mathrm{V}}_{#1}{\dsp\mathrm{[\;}}#2{\mathrm{\;]}}\fi}

\def\sg{Sierpi\'{n}ski gasket}
\def\parr{\par\noindent}
\def\eqsb{\begin{eqnarray*}}\def\eqse{\end{eqnarray*}}

\def\defd{\mathop{\stackrel{\rm def}{=}}}

\title{
 Loop-erased random walk on the Sierpinski gasket
}
\author{
Kumiko Hattori, 
Michiaki Mizuno
\\ 
\small Department of Mathematics and Information Sciences,
\\
\small  Tokyo Metropolitan University, Hachioji, Tokyo 192-0397, Japan.
\\ \small email: \url{khattori@tmu.ac.jp}
\\ \and
} 

\begin{document}
\maketitle

\thispagestyle{myheadings}
\begin{abstract}
We consider a  model of loop-erased random walks on the finite pre-{\sg} which permits rigorous 
analysis.  
We prove the existence of the scaling limit and show that the path of the 
limiting process is almost surely self-avoiding, while having Hausdorff dimension strictly greater than $1$.
This result means that the path 
has infinitely fine creases, while having no self-intersection.  
Our loop-erasing procedure is formulated 
by  a `larger-scale-loops-first' rule.  It enables us to obtain exact recursion relations, 
making use of `self-similarity' of a fractal structure.

\end{abstract}



\section{Introduction}
\seca{introduction}

In this paper, we consider a  model of loop-erased random walks on the finite pre-{\sg} which permits rigorous 
analysis.  

A loop-erased random walk is a kind of self-avoiding walk, which is a random walk that cannot visit 
any point more than once.  Concerning self-avoiding walks,  
there have been questions that are simple to ask but difficult to answer, such as: 
How far can an $n$-step self-avoiding walk go in average?  
Does it have a scaling-limit?
The non-Markov property of the walk makes the matter so difficult that we still do not know 
rigorous proofs for the `standard' model on the low-dimensional (2- and 3- dimensional) square lattices,  
which corresponds to the uniform measure on self-avoiding paths of a given length 
(\cite{ms}). 
As such, we believe a self-avoiding walk on the pre-{\sg} (a lattice version) 
serves as an interesting low-dimensional model, since it is solvable.

In \cite{hhk1, hh, hhk2, HHH}, 
models for self-avoiding walks on the 
2- and 3-dimensional pre-{\sg} were investigated, and a positive answer to the second question, above,  
was established; in addition, some path properties of the limit process 
were proved such as Hausdorff dimensions, H\"older continuity,  
whether the limit is also self-avoiding, and so on. 
In \cite{hk, hhk2}, some results were provided with regard to the first question.  The values of 
the mean-square displacement exponents obtained earlier by scaling arguments in physics literature 
were proved.  

On the other hand, Lawler \cite{lawler} defined a loop-erased random walk on square lattices, 
which is a process obtained by  chronologically 
erasing the loops from a simple random walk.  It is another kind of self-avoiding walk, 
but in this case, one can make use of the properties of simple random walks, on which there has been much study,  for analysis.
The scaling limit of the loop-erased random walk on the 2-dimensional lattice 
has been studied, using Schramm Loewner Evolution (SLE).  To name a few works in this line, \cite{LSW}, \cite{schramm}.
In \cite{Kozma}, Kozma proved the existence of the scaling limit of  the 3-dimensional loop-erased random walk.  


In this paper, we define a loop-erased random walk on the pre-{\sg} by employing 
a `larger-scale-loops-first' rule, which enables us to obtain recursion relations, 
making use of `self-similarity' 
of a fractal structure, instead of translational 
invariance of the square lattices. 
Our loop-erased walk will also be self-avoiding, but we shall show that it  belongs to a different universality class 
from the self-avoiding walk with uniform measure.
We shall also prove the existence of the scaling limit, and that the path of the 
limiting process is almost surely self-avoiding, while having Hausdorff dimension 
$\log \{ \frac{1}{15}(20+\sqrt{205})\} /\log 2$=1.1939 \ldots . This result means that the path 
has infinitely fine creases, while having no self-intersection.

Shinoda \cite{shinoda} obtained the exponent for the mean-square displacement for loop-erased random walks
on the pre-{\sg} through uniform spanning trees.  
In the physics literature, D. Dhar and A. Dhar \cite{DD} investigated 
the distribution of sizes of erased loops in terms of spanning tree and scaling arguments.  
Our path Hausdorff dimension is consistent with their results, so it is our belief that 
our larger-scale-loops-first formulation  is a natural procedure to study.

In Section 2, we describe the set-up of our model and the loop-erasing procedure, and show that 
the asymptotics of path length is consistent with the results in \cite{DD} and \cite{shinoda}.  
Section 3 is devoted to the examination of scaling limit. 

\vspace{0.5cm}\parr
{\bf Acknowledgement}
One of the authors (K. Hattori) is supported by the Grant-in-Aid for Scientific Research (C), Japan Society for the 
Promotion of Science.  We would like to thank T. Hattori and  S. Horocholyn  for helpful discussion and advice, 
and T. Itani for technical assistance.

\section{Paths on the pre-{\sg}s}
\seca{2}
\subsection{The pre-{\sg}s.}


We consider  the pre-{\sg}, a lattice version of the {\sg}, which is a fractal with 
Hausdorff dimension $\log 3/\log 2$.  (For fractals, see \cite{falconer}.)
Let us recall the definition of the pre-{\sg}: by denoting  
$O=(0,0),\,
a_0=({1 \over 2},{\sqrt{3} \over 2}),\,
b_0=(1,0)\, $, 
and for each $N \in 
{\mathbb N}$, 
$a_N=2^Na_0, \ \ \ b_N=2^Nb_0,$  
then define $F'_{0}$ be the graph that consists of three vertices and three edges of $\triangle Oa_0b_0 $
and define the recursive sequence of graphs 
$\{F'_{N}\}_{N=0}^{\infty }$
by
\\
\unitlength 0.1in
\begin{picture}( 42.8000, 20.2400)(  1.9500,-22.7200)
%
{\color[named]{Black}{%
\special{pn 8}%
\special{pa 2168 2042}%
\special{pa 1218 2042}%
\special{pa 1694 1220}%
\special{pa 2168 2042}%
\special{pa 1218 2042}%
\special{fp}%
}}%
%
{\color[named]{Black}{%
\special{pn 8}%
\special{pa 242 2038}%
\special{pa 480 1628}%
\special{pa 718 2038}%
\special{pa 242 2038}%
\special{pa 480 1628}%
\special{fp}%
}}%
%
{\color[named]{Black}{%
\special{pn 8}%
\special{pa 1694 2042}%
\special{pa 1932 1630}%
\special{pa 2170 2042}%
\special{pa 1694 2042}%
\special{pa 1932 1630}%
\special{fp}%
}}%
%
{\color[named]{Black}{%
\special{pn 8}%
\special{pa 1694 2042}%
\special{pa 1932 1630}%
\special{pa 2170 2042}%
\special{pa 1694 2042}%
\special{pa 1932 1630}%
\special{fp}%
}}%
%
{\color[named]{Black}{%
\special{pn 8}%
\special{pa 1694 2042}%
\special{pa 1932 1630}%
\special{pa 2170 2042}%
\special{pa 1694 2042}%
\special{pa 1932 1630}%
\special{fp}%
}}%
%
{\color[named]{Black}{%
\special{pn 8}%
\special{pa 1694 2042}%
\special{pa 1932 1630}%
\special{pa 2170 2042}%
\special{pa 1694 2042}%
\special{pa 1932 1630}%
\special{fp}%
}}%
%
{\color[named]{Black}{%
\special{pn 8}%
\special{pa 1218 2042}%
\special{pa 1456 1630}%
\special{pa 1694 2042}%
\special{pa 1218 2042}%
\special{pa 1456 1630}%
\special{fp}%
}}%
%
{\color[named]{Black}{%
\special{pn 8}%
\special{pa 1218 2042}%
\special{pa 1456 1630}%
\special{pa 1694 2042}%
\special{pa 1218 2042}%
\special{pa 1456 1630}%
\special{fp}%
}}%
%
{\color[named]{Black}{%
\special{pn 8}%
\special{pa 1218 2042}%
\special{pa 1456 1630}%
\special{pa 1694 2042}%
\special{pa 1218 2042}%
\special{pa 1456 1630}%
\special{fp}%
}}%
%
{\color[named]{Black}{%
\special{pn 8}%
\special{pa 1218 2042}%
\special{pa 1456 1630}%
\special{pa 1694 2042}%
\special{pa 1218 2042}%
\special{pa 1456 1630}%
\special{fp}%
}}%
%
{\color[named]{Black}{%
\special{pn 8}%
\special{pa 1456 1630}%
\special{pa 1694 1218}%
\special{pa 1932 1630}%
\special{pa 1456 1630}%
\special{pa 1694 1218}%
\special{fp}%
}}%
\put(21.5600,-20.8800){\makebox(0,0)[lt]{$b_1$}}%
\put(14.5200,-16.2900){\makebox(0,0)[rb]{$a_0$}}%
\put(16.5000,-12.1500){\makebox(0,0)[lb]{$a_1$}}%
\put(16.9300,-21.4500){\makebox(0,0){$b_0$}}%
\put(11.7100,-20.8800){\makebox(0,0)[lt]{$O$}}%
\put(1.9500,-20.8500){\makebox(0,0)[lt]{$O$}}%
\put(7.4300,-21.4700){\makebox(0,0){$b_0$}}%
\put(5.5700,-16.1400){\makebox(0,0)[rb]{$a_0$}}%
\put(16.9300,-23.3700){\makebox(0,0){$F'_1$}}%
%
{\color[named]{Black}{%
\special{pn 8}%
\special{pa 3526 2042}%
\special{pa 2576 2042}%
\special{pa 3050 1220}%
\special{pa 3526 2042}%
\special{pa 2576 2042}%
\special{fp}%
}}%
%
{\color[named]{Black}{%
\special{pn 8}%
\special{pa 3050 2042}%
\special{pa 3288 1630}%
\special{pa 3526 2042}%
\special{pa 3050 2042}%
\special{pa 3288 1630}%
\special{fp}%
}}%
%
{\color[named]{Black}{%
\special{pn 8}%
\special{pa 3050 2042}%
\special{pa 3288 1630}%
\special{pa 3526 2042}%
\special{pa 3050 2042}%
\special{pa 3288 1630}%
\special{fp}%
}}%
%
{\color[named]{Black}{%
\special{pn 8}%
\special{pa 3050 2042}%
\special{pa 3288 1630}%
\special{pa 3526 2042}%
\special{pa 3050 2042}%
\special{pa 3288 1630}%
\special{fp}%
}}%
%
{\color[named]{Black}{%
\special{pn 8}%
\special{pa 3050 2042}%
\special{pa 3288 1630}%
\special{pa 3526 2042}%
\special{pa 3050 2042}%
\special{pa 3288 1630}%
\special{fp}%
}}%
%
{\color[named]{Black}{%
\special{pn 8}%
\special{pa 2576 2042}%
\special{pa 2814 1630}%
\special{pa 3052 2042}%
\special{pa 2576 2042}%
\special{pa 2814 1630}%
\special{fp}%
}}%
%
{\color[named]{Black}{%
\special{pn 8}%
\special{pa 2576 2042}%
\special{pa 2814 1630}%
\special{pa 3052 2042}%
\special{pa 2576 2042}%
\special{pa 2814 1630}%
\special{fp}%
}}%
%
{\color[named]{Black}{%
\special{pn 8}%
\special{pa 2576 2042}%
\special{pa 2814 1630}%
\special{pa 3052 2042}%
\special{pa 2576 2042}%
\special{pa 2814 1630}%
\special{fp}%
}}%
%
{\color[named]{Black}{%
\special{pn 8}%
\special{pa 2576 2042}%
\special{pa 2814 1630}%
\special{pa 3052 2042}%
\special{pa 2576 2042}%
\special{pa 2814 1630}%
\special{fp}%
}}%
%
{\color[named]{Black}{%
\special{pn 8}%
\special{pa 2814 1630}%
\special{pa 3050 1218}%
\special{pa 3288 1630}%
\special{pa 2814 1630}%
\special{pa 3050 1218}%
\special{fp}%
}}%
\put(34.9400,-20.8800){\makebox(0,0)[lt]{$b_1$}}%
\put(28.0900,-16.2900){\makebox(0,0)[rb]{$a_0$}}%
\put(29.2600,-12.1000){\makebox(0,0)[lb]{$a_1$}}%
\put(30.6500,-21.4500){\makebox(0,0){$b_0$}}%
\put(25.2800,-20.8800){\makebox(0,0)[lt]{$O$}}%
\put(35.2800,-23.3700){\makebox(0,0){$F'_2$}}%
%
{\color[named]{Black}{%
\special{pn 8}%
\special{pa 4476 2040}%
\special{pa 3526 2040}%
\special{pa 4000 1216}%
\special{pa 4476 2040}%
\special{pa 3526 2040}%
\special{fp}%
}}%
%
{\color[named]{Black}{%
\special{pn 8}%
\special{pa 4000 2040}%
\special{pa 4238 1628}%
\special{pa 4476 2040}%
\special{pa 4000 2040}%
\special{pa 4238 1628}%
\special{fp}%
}}%
%
{\color[named]{Black}{%
\special{pn 8}%
\special{pa 4000 2040}%
\special{pa 4238 1628}%
\special{pa 4476 2040}%
\special{pa 4000 2040}%
\special{pa 4238 1628}%
\special{fp}%
}}%
%
{\color[named]{Black}{%
\special{pn 8}%
\special{pa 4000 2040}%
\special{pa 4238 1628}%
\special{pa 4476 2040}%
\special{pa 4000 2040}%
\special{pa 4238 1628}%
\special{fp}%
}}%
%
{\color[named]{Black}{%
\special{pn 8}%
\special{pa 4000 2040}%
\special{pa 4238 1628}%
\special{pa 4476 2040}%
\special{pa 4000 2040}%
\special{pa 4238 1628}%
\special{fp}%
}}%
%
{\color[named]{Black}{%
\special{pn 8}%
\special{pa 3526 2040}%
\special{pa 3762 1628}%
\special{pa 4000 2040}%
\special{pa 3526 2040}%
\special{pa 3762 1628}%
\special{fp}%
}}%
%
{\color[named]{Black}{%
\special{pn 8}%
\special{pa 3526 2040}%
\special{pa 3762 1628}%
\special{pa 4000 2040}%
\special{pa 3526 2040}%
\special{pa 3762 1628}%
\special{fp}%
}}%
%
{\color[named]{Black}{%
\special{pn 8}%
\special{pa 3526 2040}%
\special{pa 3762 1628}%
\special{pa 4000 2040}%
\special{pa 3526 2040}%
\special{pa 3762 1628}%
\special{fp}%
}}%
%
{\color[named]{Black}{%
\special{pn 8}%
\special{pa 3526 2040}%
\special{pa 3762 1628}%
\special{pa 4000 2040}%
\special{pa 3526 2040}%
\special{pa 3762 1628}%
\special{fp}%
}}%
%
{\color[named]{Black}{%
\special{pn 8}%
\special{pa 3762 1628}%
\special{pa 4000 1216}%
\special{pa 4238 1628}%
\special{pa 3762 1628}%
\special{pa 4000 1216}%
\special{fp}%
}}%
\put(44.6300,-20.8700){\makebox(0,0)[lt]{$b_2$}}%
%
{\color[named]{Black}{%
\special{pn 8}%
\special{pa 4000 1216}%
\special{pa 3050 1216}%
\special{pa 3526 394}%
\special{pa 4000 1216}%
\special{pa 3050 1216}%
\special{fp}%
}}%
%
{\color[named]{Black}{%
\special{pn 8}%
\special{pa 3050 1216}%
\special{pa 3288 804}%
\special{pa 3526 1216}%
\special{pa 3050 1216}%
\special{pa 3288 804}%
\special{fp}%
}}%
%
{\color[named]{Black}{%
\special{pn 8}%
\special{pa 3524 1216}%
\special{pa 3762 804}%
\special{pa 4000 1216}%
\special{pa 3524 1216}%
\special{pa 3762 804}%
\special{fp}%
}}%
%
{\color[named]{Black}{%
\special{pn 8}%
\special{pa 3286 804}%
\special{pa 3526 390}%
\special{pa 3762 804}%
\special{pa 3286 804}%
\special{pa 3526 390}%
\special{fp}%
}}%
\put(34.4300,-3.7800){\makebox(0,0)[lb]{$a_2$}}%
\put(4.7600,-23.3600){\makebox(0,0){$F'_0$}}%
\end{picture}%
 \\
\\
\hspace{3cm} {\bf Fig. 1}
\vspace{0.5cm}
\[
F'_{N+1}=F'_{N} \cup (F'_{N}+a_N) \cup (F'_{N}+b_N)
, \qquad N \in {\mathbb Z}_+ =\{0,1,2, \ldots \}\, ,
\]
where $A+a=\{x+a\ :\  x\in A\}$ and $kA=\{kx\ :\  x\in A\}$.  
$F'_0$, $F'_1$ and $F'_2$ are shown in Fig. 1. 

Finally, we let $F''_N$ be the union of $F'_N$ and its reflection with respect to the $y$-axis, 
and denote $\dsp F_{0} =\bigcup _{N=1}^{\infty} F''_N$;
the graph $F_{0}$ is called the (infinite) {\bf pre-{\sg}}.  $F_0$ is shown in Fig. 2.
\\
\input{Fig2.tex}\\
\\
\hspace{3cm} {\bf Fig. 2}
\vspace{0.5cm}

Furthermore, by letting 
$G_0$ and $E_0$ denote the set of vertices and the set of edges of $F_0$, respectively,
we see that, for each $N \in {\mathbb Z}_+$, $F_{N}=2^N F_{0}$ 
can be regarded as a coarse graph with vertices $G_N=\{2^Nx\ :\ x\in G_0 \}$
and edges $E_N=\{ 2^N\overline{xy} \ :\ \overline{xy} \in E_0\}$.
Given $x\in G_N$, let ${\cal N}_N(x)$ be the four nearest neighbors of 
$x$ on $F_N$, that is, ${\cal N}_N(x)=\{y \in G_N :\ \overline{xy} \in E_N\}$.

\vspace{0.5cm}\parr
\subsection{Paths on the pre-{\sg}s.}

Let us denote the set of finite paths on $F_0$ by
\[W = \{\ w=(w(0), w(1), \cdots , w(n)):
\  w(0) \in G_0, \ w(i)\in {\cal N}_{0}(w(i-1)), 
\ 1 \leq i \leq n,  \ n \in {\mathbb N}  \},\]
and the set of finite paths on $F_0$ starting at $O$ by
\[W^*=\{\  w\in W\ :\ w(0)=O\ \}.\]
This gives the natural definition for 
the length $\ell $ of a path  $w=(w(0), w(1), \cdots , w(n))\in W$; namely, $\ell (w)=n$. 
 
For a path $w\in W$ and $A \subset G_0$,
 we define the  hitting time of $A$ 
 by
\[T_A(w)=\inf  \{j \geq 0 :\ w(j) \in A\},\]
where we set $\inf \emptyset =\infty $. 
By taking 
$w\in W$ and $M \in \pintegers$, we shall define the recursive sequence 
$\{T_i^M(w)\}_{i=0}^{m}$ of
{\bf hitting times of $G_M$} as follows:
Let 
 $T_{0}^{M}(w)=T_{G_M}$, and for $i\geq 1$, let
\[T_{i}^{M}(w)=\inf \{j>T_{i-1}^{M}(w) :\  w(j)\in G_{M}\setminus
\{w(T_{i-1}^{M}(w))\}\};\]
here we take  $m$ to be the smallest integer such that 
$T_{m+1}^{M}(w)=\infty $. Then 
$T_{i}^{M}(w)$ can be interpreted as being the time taken for the path $w$ to 
hit vertices in $G_{M}$ for the $(i+1)$-th time, 
under the condition that if $w$ hits the same vertex in $G_{M}$ 
more than once in a row, we count it only once.

Now we consider two sequences of subsets of $W^*$ as follows: 
for each $N \in \pintegers$, let the set of paths from $O$ to $a_N$, which  
do not hit any other vertices in $G_N$ on the way, be 
\[W_N=\{ w =(w(0),w(1),\cdots,w(n)) \in W^* : \  w(n)=a_N, \ n=T_1^N(w) \},\]  
and let the set of paths from from $O$ to $a_N$ that 
hit $b_N$ `once' on the way (subject to the counting rule explained above) be 
\[V_N=\{w= (w(0),w(1),\cdots,w(n))\in W^* : \  w(n)=a_N,  w(T_1^N(w))=b_N, 
n=T_2^N(w)\}.\]

Then for a path $w\in W$ and $M\in \pintegers$, we 
 define the  {\bf coarse-graining map}
$Q_{M}$ by \[(Q_{M}w)(i)=w(T_{i}^{M}(w)), \ \ \mbox{ for } i=0,1,2,\ldots, m,\]
where $m$ is the smallest integer such that 
$T_{m+1}^{M}(w)=\infty $. Thus,
\[Q_M w=[w(T_0^M(w)), w(T_1^M(w)), \ldots , 
w(T_m^M(w)) ]\]
is  a path on a coarser graph $F_M$. For 
$w\in W_N\cup V_N$ and $M \leq N$, the end point of the coarse-grained path is $w(T_m^M(w))=a_N$, 
and if we write $(2^{-M}Q_{M}w)(i)=2^{-M}w(T_{i}^{M}(w))$,
then $2^{-M}Q_{M}w$ is a path in $W_{N-M}\cup V_{N-M}$ and $\ell (2^{-M}Q_{M}w)=m$.
Notice that if $ M \leq N$, then $Q_{N}\circ Q_{M}=Q_{N}$.
Throughout the following, we write simply $w(T^M_i)$ instead of $w(T^M_i(w))$.

\subsection{Loop-erased paths.}
\seca{2.3}

Let $\Gamma $ be the set of self-avoiding paths starting at $O$:
\[\Gamma = \{\ (w(0), w(1), \cdots , w(n)) \in W^*:
\  w(i) \neq w(j), \ i \neq j,\ n\in \nintegers \  \} ,
\]
and let us denote the following  two subsets of $\Gamma $ :
\[\hat{W}_N =W_N\cap \Gamma ,\ \  \hat{V}_N=V_N \cap \Gamma .\]

For  $(w(0), w(1), \cdots , w(n))\in W^*$, 
We call a path segment $[w(i),w(i+1), \ldots , w(j)]$ a loop if 
there are $i, j$, 
$0\leq i < j \leq n$ such that $w(i)=w(j)$ and 
$w(k)\neq w(i)$ for any $i<k<j$. .  

We shall now describe a loop-erasing procedure for paths in  $W_1\cup V_1$:
\itmb
\item[(i)] Erase all the loops formed at $O$;
\item[(ii)] Progress one step forward along the path, and 
erase all the loops at the new position;
\item[(iii)] Iterate this process, taking another step forward along the path and erasing the 
loops there, until reaching $a_1$ (the endpoint of all paths in $W_1$ and $V_1$).
\itme

To be precise, for $w \in W_1\cup V_1$, define the recursive sequence $\{s_i\}_{i=0}^{m}$ 
\[s_0=\sup \{j:w(j)=O \},\]
\[s_i=\sup \{j:w(j)=w(s_{i-1}+1)\}.\]
If $s_i>s_{i-1}+1$, then 
$[w(s_{i-1}+1), w(s_{i-1}+2), \ldots , w(s_{i}-1), w(s_{i})]$
forms a loop, starting and ending at 
$w(s_{i-1}+1)=w(s_{i})$.  We erase it by removing all of the points 
$w(s_{i-1}+1), w(s_{i-1}+2), \ldots , w(s_{i}-2), w(s_{i}-1)$.
If $w(s_m)=a_1$,  then 
we have obtained a loop-erased path,
\[L w =[w(s_0), w(s_1), \ldots ,w(s_m)]\in 
\hat{W}_1 \cup \hat{V}_1.\]
Note that $w \in W_1$ implies $Lw \in \hat{W}_1$, but that 
 $w \in V_1$ can result in  $Lw \in \hat{W}_1$, with $b_1$ being erased 
 together with a loop.
So far, our loop-erasing procedure is the same as that defined for paths on 
${\mathbb Z}^d $ in \cite{lawler}. 

We shall generalize the above procedure to a loop-erasing procedure for a path $w $ in 
$W_N \cup V_N$ that yields a self-avoiding 
path in $\hat{W}_N \cup \hat{V}_N$. 
The idea is to first erase loops of `largest scale', and then go down to
`smaller scales'  step by step.  For this purpose, we need the notion of `skeletons'.

Let ${\cal T}_M$ be the set of all upward (closed and filled) triangles which are translations of 
 $\triangle Oa_Mb_M$ and whose vertices are in $G_M$; an 
element of ${\cal T}_M$ is called a {\bf $2^M$-triangle}. 
For $w \in  W$ and $M\geq 0$, we shall define 
a sequence $(\Delta_1, \ldots , \Delta_k)$ of $2^M$-triangles $w$ `passes through' and 
a sequence $\{T_i^{ex,M}(w)\}_{i=1}^{k}$ of 
exit times from them as a subsequence of  $\{T_i^{M}(w)\}_{i=1}^{m}$, as follows:
We start by defining  $T_0^{ex, M}(w)=T_0^{M}(w)$.  (Thus If $w\in W^*$, then $T_0^{ex, M}(w)=0$.)
There is a unique element of ${\cal T}_M$ that contains $w(T_0^M)$ and $w(T_1^M)$, 
which we denote by  $\Delta _1$. 
For $i\geq 1$, define 
\[j(i)=\min \{j \geq 0\ : \ j<m,\ T_j^M(w)>T_{i-1}^{ex, M}(w),\ 
w(T_{j+1}^M(w))\not\in \Delta _i\},\] 
if the minimum exists, 
otherwise $j(i)=m$.   
Then define $ T_i^{ex, M}(w)=T_{j(i)}^{M}(w)$, and let $\Delta _{i+1}$ be
the unique $2^M$-triangle that contains both $w(T_i^{ex, M})$ and $w(T_{j(i)+1}^{M})$. 
By definition, we see that  $\Delta _{i} \cap 
\Delta _{i+1 }$ is a one-point set $\{w(T_i^{ex, M})\}$, for $i=1, \ldots , k-1$. 
We denote the sequence of these triangles by 
$\sigma _M(w)=( \Delta _1, \ldots , \Delta _k)$, 
and call it the {\bf $2^M$-skeleton} of $w$.
We call the sequence  
$\{T_i^{ex, M}(w)\}_{i=0, 1, \ldots , k}$ 
{\bf exit times} from the triangles in the skeleton.  
For each $i$, there is an $n=n(i)$ such that 
$T_{i-1}^{ex, M}(w)=T^M_{n}(w)$.  
We say $\Delta _i \in \sigma _M(w)$ is an element of {\bf Type 1} if $T_{i}^{ex, M}(w)=T^M_{n+1}$,
and an element of {\bf Type 2} if  $T_{i}^{ex, M}(w)=T^M_{n+2}$. 
If $w\in \hat{W}_N\cup \hat{V}_N$ for some $N$, then $\Delta _1, \ldots , \Delta _{k}$ are mutually distinct, 
and each of them is either of Type 1 or of Type 2. 

Assume $w\in W_N\cup V_N$ for some $N$ and $M\leq N$.
For each $\Delta $ in $\sigma _M(w)$, the {\bf path segment of $w$ in $\Delta $} is 
\[[w(n), \ T^{ex, M}_{i-1}(w) \leq n \leq T^{ex, M}_{i}],\] 
and it is denoted by $w|_{\Delta }$. 
Note that the definition of $T^M_i$'s allows a path segment $w|_{\Delta }$ to leak into two neighboring 
$2^M$-triangles.  
It should be noted that the subgraph contained in $\Delta $ and its neighboring triangles has the same structure as 
$\triangle Oa_Mb_M$ and its neighbors, which implies that $w|_{\Delta }$ can be naturally identified 
with some path in $\triangle Oa_Mb_M$ and its neighbors starting at $O$, 
by translation, rotation and reflection.  For convenience we 
shall denote this identification by $\eta $, and write: 
\eqnb \eqna{eta} \eta (w|_{\Delta }) =v \in W_M \cup V_M, \eqne
where the entrance to $\Delta $ is mapped to $O$ and the exit to $a_M$. 



To introduce the loop-erasing operation for paths in $W_N\cup V_N$, let us take 
a loop $[w(i), w(i+1), \ldots , w(i+i_0)]$ that is contained in  $w \in W_N\cup V_N$
, and  define its diameter by $d=\max \{i<j\leq i+i_0:|w(j)-w(i)|\}$.
The loop $[w(i), w(i+1), \ldots , w(i+i_0)]$ is said to be a {\bf $2^{M}$-scale loop}, 
whenever  there exists an $M\in {\mathbb Z}_+$ such that 
\[\max \{N':w(i)=w(i+i_0)\in G_{N'}\}=M \mbox{ and } d\geq 2^M.\]
Then the definition implies that  $w $ has a $2^{N-1}$-scale 
loop if and only if the coarse-grained path $Q_{N-1}w$ has a loop. 
The operation of erasing largest-scale loops can be reduced to erasing loops from a path in $W_1\cup V_1$, 
which we shall show below by induction.


Let $w \in W_N \cup V_N$ (Fig. 3(a)). 
we define the operation of `erasing the largest-scale loops'  as follows:
\itmb
\item[1)] Coarse-grain $w$ to obtain 
\[w'=Q_{N-1}w =[w(T^{N-1}_{0}), w(T^{N-1}_{1}), \ldots , 
w(T^{N-1}_{k})],\]
where $w(T_k^{N-1})=a_N$ (Fig. 3(b)).  
We note that $2^{-(N-1)}w' \in W_{1} \cup V_{1}$.

\item[2)] Similarly to the procedure for $W_1\cup V_1$, 
erase loops from $w'$, using the following sequence and defining the mapping $L$: 
\[s_0=\sup \{j:w(T_j^{N-1})=O \},\]
\[s_i=\sup \{j:w(T_j^{N-1})=w(T_{s_{i-1}+1}^{N-1})\},\ i \geq 1, \]
and 
\[Lw'=[w(T^{N-1}_{s_0}), w(T^{N-1}_{s_1}), \ldots , 
w(T^{N-1}_{s_m})],\]
where $w(T^{N-1}_{s_m})=a_N$ (Fig. 3(c)).
We note here that $2^{-(N-1)}Lw' \in \hat{W}_1\cup \hat{V}_1$. 
 
\item[3)] Make a path by concatenation of $m$ parts chosen from the 
original path ; 
\[L_{N-1}w=[w _0, w _1, \ldots , w_{m-1}, a_N],\]
where 
\[w_i=[w(T_{s_i}^{N-1}), w (T_{s_i}^{N-1}+1)
\ldots , w(T_{s_i+1}^{N-1}-1)],\ \ 
i=0, \cdots , m-1.\]
\itme

By steps  1)--3), we have obtained $L_{N-1}w \in W_N\cup V_N$ with all $2^{N-1}$-scale loops of $w$ erased
(Fig. 3(d)).
\\
\\
\input{Fig3.tex}\\
\\
\hspace{3cm} {\bf Fig. 3}
\vspace{0.5cm}

Using above as a base step, we shall now describe the induction step of our operation: 
Let $w\in W_N \cup V_N$. 
For $M\leq N$,  
assume that all of the $2^N$- to $2^{M}$-scale loops have been erased from the 
path $w$, and denote  the resulting path  $w'$, and its $2^M$-skeleton by $\sigma _M(w')$.
Additionally, for  each $\Delta \in \sigma _M(w')$, we shall (implicitly) use the  identification $\eta $ defined in 
\eqnu{eta} to identify  
$Q_{M-1}w'|_{\Delta}$ with a path in $W_1\cup V_1$. 
 
\itmb
\item[L1)] 
Coarse-grain $w'$ to obtain $Q_{M-1}w'$ and consider 
\[Q_{M-1}w' |_{\Delta }=[w'(T^{M-1}_{k}), w'(T^{M-1}_{k+1}), \ldots , 
w'(T^{M-1}_{k+k_0})], \] 
where $w'(T^{M-1}_{k})$ is the entrance point to $\Delta $ and $w'(T^{M-1}_{k+k_0})$
the exit point from $\Delta $. 

\item[L2)] Erase loops from $Q_{M-1}w' |_{\Delta }$ as in the procedure 
for $W_1\cup V_1$ by defining the sequence $\{s_i\}_{i=1}^m$ by 
\[s_0=\sup \{j:w'(T_j^{M-1})=w'(T_k^{M-1}) \},\]
\[s_i=\sup \{j:w'(T_j^{M-1})=w'(T_{s_{i-1}+1}^{M-1})\},\ i \geq 1, \]
and denoting 
\[L(Q_{M-1}w'|_{\Delta })=[w'(T^{M-1}_{s_0}), w'(T^{M-1}_{s_1}), \ldots , 
w'(T^{M-1}_{s_m})],\]
where $w'(T^{M-1}_{s_0})= w'(T^{M-1}_{k})$ and $w'(T^{M-1}_{s_m})=w'(T^{M-1}_{k+k_0})$.

\item[L3)] Make a path segment in $\Delta $ by concatenation of $m$ parts chosen from the 
original path and the exit point and denote it by 
\[L_{M-1}(w|_{\Delta })=[w' _0, w' _1, \ldots , w'_{m-1}, w'(T_{s_m}^{M-1})],\]
where 
\[w'_i=[w'(T_{s_i}^{M-1}), w' (T_{s_i}^{M-1}+1)
\ldots , w'(T_{s_i+1}^{M-1}-1)],\ \ 
i=0, \cdots , m-1.\]

\item[L4)]
Make a whole path $w''=L_{M-1}w$ by concatenation of parts obtained in L3) over all  $\Delta \in \sigma _M(w')$.
\itme

Thus, by the procedure above, we have erased all of the $2^{M-1}$-scale loops from 
$w$.
Now denote by $\hat{Q}_{M-1}w $ the path obtained by concatenation of $L(Q_{M-1}w'|_{\Delta})$ 
obtained in L2); then   
it is a path on $F_{M-1}$, in the sense that $Q_{M-1}(\hat{Q}_{M-1}w)=\hat{Q}_{M-1}w $, 
from $O$ to $a_N$ without loops.
Observe that $\hat{Q}_{M-1}w=Q_{M-1}w''$.  Although  it may occur that 
$\sigma _{M-1}(w'') \neq \sigma _{M-1}(w')$,  it holds that $\sigma _M(w'')=\sigma _M(w')$, which
can be extended to $\sigma _K(w')=\sigma _K(w'')$ for any $K\geq M$.

We then continue this operation until we have erased all of the loops and have $Lw=L_0w=\hat{Q}_0w$.
Thus, by construction, our loop-erasing operation is essentially a repetition of loop-erasing  for 
$W_1\cup V_1$.
We remark that  the procedure implies that for any  $w\in W_N \cup V_N$,

\eqnb
\eqna{sigmainv}  
\sigma _K(\hat{Q}_{M}w)=\sigma _K (\hat{Q}_{K}w) \ \mbox{ for any } M \leq K \leq N.\eqne
i.e., in the process of loop-erasing, once loops of $2^K$-scale and greater 
have been erased, the $2^K$-skeleton does not change any more.  However it should be noted 
that the types of the triangles can change from Type 2 to Type 1.


\vspace{0.5cm}\par\noindent
\subsection{Loop-erased random walks on the pre-{\sg}s.}


Let $(\tilde{\Omega }, {\cal F}, P)$ be a probability space.
A simple random walk on $F_0$ is a $G_0$-valued Markov chain 
$\{Z(i) : i\in \pintegers \}$ with transition probabilities 
\[P[Z(i+1)=y\ |\ Z(i)=x]=
\left\{
\begin{array}{ll}
\frac{1}{4} &\mbox{   if } y \in {\cal N}_0(x)
\\
0 & \mbox{   otherwise.}
\end{array}
\right.
\]
Throughout this paper, we will consider random walks starting at O, so 
finite random walk paths are elements of $W^*$, and thus, $T_i^N$'s and $Q_NZ$ can be
defined.

Consider two kinds of random walks stopped at $a_N$: one conditioned on  $Z(T_1^N)=a_N$ 
(before hitting other $G_N$ vertices), called $X_N$, and the other conditioned on 
$Z(T_1^N)=b_N$ and $Z(T_2^N)=a_N$, i.e. hitting $b_N$ on the way to $a_N$, called $X'_N$.   
These random walks then induce measures $P_N$ and $P'_N$ on $W^*$  
with support on $W_N$ and $V_N$, respectively, namely, 

\begin{eqnarray*}
P_N[w]&=&
P[X_N(i)=w(i), i=0, 1, \ldots , \ell (w)]\\
&=&
P[Z(i)=w(i), i=0, 1, \ldots , \ell (w) | Z(T_1^N)=a_N ], \ w\in W_N,
\end{eqnarray*}

\begin{eqnarray*}
P'_N[w]
&=&
P[X'_N(i)=w(i), i=0, 1, \ldots , \ell (w)]\\
&=&
 P[Z(i)=w(i), i=0, 1, \ldots , \ell (w) | Z(T_1^N)=b_N, \ \ Z(T_2^N)=a_N ],\ \ 
w \in V_N.
\end{eqnarray*}
Note that by symmetry:
\[P[ Z(T_1^N)=a_N]=1/4, \ \ \ P[Z(T_1^N)=b_N, \ \ Z(T_2^N)=a_N]=1/16.\]

Throughout this paper,  
the following propositions on the simple random walks on the pre-{\sg} will be used;   
They are straightforward consequences of  the `self-similarity', that is, $2^{-M}F_M=F_0$, 
and the property that if $x_0\in G_M$ for some $M \in {\mathbb Z}_+$, then 
for each $x\in {\cal N}_M(x_0)$ 
\[
P[Z(T_{i+1}^M)=x| Z(T_i^M)=x_0 ]=\frac{1}{4}
\]
holds. 
(For details of random walks on the {\sg}, we refer to  \cite{BP}.)



\prpb
\prpa{q1}
If $M\leq N$, then the distributions of $2^{-M}Q_MX_N$ and $2^{-M}Q_MX'_N$ 
are equal to $P_{N-M}$ and $P'_{N-M}$, respectively; in other words, $Q_MX_N$ and $Q_MX'_N$ are simple random 
walks on a coarse graph $F_M$ stopped at $a_N$. 
\prpe

Let $\eta $ be the identification map defined in the last subsection. 


\prpb
\prpa{q2}
Let $M\leq N$, and consider random walk segments conditioned on $Q_MX_N$ between the hitting times, 
\[Z_i=[X_N(t), \ T_i^M(X_N)\leq t \leq  T_{i+1}^M(X_N) ],\ \  \ i=1, \ldots , m ,\]
where $X_N( T_m^M)=a_N$.  
Then $Z_i$, $i=1, \ldots , m$ , when identified with paths in 
$W_{N-M}$ by appropriate translation, rotation and reflection, 
are independent and have the same distribution as $X_{N-M}$. 
\prpe

By applying loop-erasing operation to random walks  $X_N$ and  $X'_N$, we 
induce measures $\hat{P}_N=P_N\circ L^{-1}$ supported on $\hat{W}_N$, 
and $\hat{P}'_N=P'_N\circ L^{-1}$ supported  
on $\hat{W}_N\cup \hat{V}_N$, respectively.
Paths in $\hat{W}_1$ and $\hat{V}_1$ are  shown in Fig. 4.
\\
\input{Fig4.tex} \\
\\

\hspace{3cm} {\bf Fig. 4}
\vspace{0.5cm}

Their probabilities under  
$\hat{P}_1$ and $\hat{P}'_1$, respectively, can be obtained by direct calculation:

\[\hat{P}_1[w_1]=\frac{1}{2}, \hat{P}_1[w_2]=\frac{2}{15}, \hat{P}_1[w_3]=\frac{2}{15}, 
\hat{P}_1[w_4]=\frac{1}{30}, \hat{P}_1[w_5]=\frac{1}{30}, \hat{P}_1[w_6]=\frac{1}{30}, 
\hat{P}_1[w_7]=\frac{2}{15}, \]
\[\hat{P}'_1[w_1]=\frac{1}{9}, \hat{P}'_1[w_2]=\frac{11}{90}, \hat{P}'_1[w_3]=\frac{11}{90}, 
\hat{P}'_1[w_4]=\frac{2}{45}, \hat{P}'_1[w_5]=\frac{2}{45}, \hat{P}'_1[w_6]=\frac{2}{45},\] 
\[\hat{P}'_1[w_7]=\frac{8}{45}, \hat{P}'_1[w_8]=\frac{2}{9}, \hat{P}'_1[w_9]=\frac{1}{18}, 
\hat{P}'_1[w_{10}]=\frac{1}{18}.\]

\vspace{0.5cm}\parr

For $w\in \hat{W}_N\cup \hat{V}_N$, let us denote  the number of Type 1 triangles and
Type 2  triangles in $\sigma _0(w)$  by $s_1(w)$ and $s_2(w)$, respectively.
(This implies that $\ell (w)=s_1 (w)+2 s_2(w)$.) 
Define two sequences, $\{\Phi _N\}_{N \in {\mathbb N}}$ and $\{\Theta _N\}_{N \in {\mathbb N}}$, of 
generating functions by:  
\[\Phi _N (x,y)=\sum _{w\in \hat{W}_N }\hat{P}_N(w)x^{s_1(w)}y^{s_2(w)},\]
\[\Theta _N (x,y)=\sum _{w\in \hat{V}_N }\hat{P}'_N(w)x^{s_1(w)}y^{s_2(w)},\ \ \ x,y\geq 0.\]
For simplicity, we shall denote $\Phi _1 (x,y)$ and $\Theta _1 (x,y)$ by $\Phi (x,y)$ and $\Theta (x,y)$.



\prpb
\prpa{recursion}
The above generationg functions satisfy the following recursion relations for all $N \in {\mathbb N}$:
\[\Phi (x,y)=\frac{1}{30}(15x^2+8xy+y^2+2x^2y+4x^3).
\]
\[\Theta (x,y)=\frac{1}{45}(5x^2+11xy+2y^2+14x^2y+8x^3+5xy^2).
\]
\[\Phi _{N+1}(x,y)=\Phi _{N}(\Phi (x,y), \Theta (x,y)).\]
\[\Theta _{N+1}(x,y)=\Theta _{N}(\Phi (x,y), \Theta (x,y)).\]
\prpe

\prfb

We shall first express $\hat{P}_{N+1}$ in terms of $\hat{P}_{N}$, $\hat{P}_{1}$ and $\hat{P}'_{1}$.
If we recall the procedure for  obtaining $\hat{Q}_1X_{N+1}$ from $X_{N+1}$, 
we notice that it is the same as the procedure to obtain $LX_N$ from $X_N$, except 
that everything is twice larger in the case of $X_{N+1}$. 
This together with  
\prpu{q1} implies that the distribution of  $2^{-1}\hat{Q}_1X_{N+1}$ is equal to  $\hat{P}_{N}$, 
namely,  
\[P_{N+1}[\ v \ :\  \frac{1}{2}\hat{Q}_1v=u\ ] =\hat{P}_N[\ u\ ].\]

On the other hand, we have from \eqnu{sigmainv} 
\[\sigma _1(\hat{Q}_1X_{N+1})=\sigma _1(LX_{N+1}).\]
The rest of the loop-erasing procedure to obtain $LX_{N+1}$ together with 
\prpu{q2} implies that conditioned on $\hat{Q}_1X_{N+1}$, 
the walk segments of  $L_1X_{N+1}$ in $\Delta  \in \sigma _1(\hat{Q}_1X_{N+1})$ 
have the same distribution as either $X_1$ or $X'_1$ (modulo appropriate transformation), 
and that they are mutually independent, 
which further implies that 
$LX_{N+1}|_{\Delta }$ are independent. 

Keeping these observations in mind, we calculate $\hat{P}_{N+1}[w]$ for 
$w\in \hat{W}_{N+1}$. Let  
$\sigma _1(w)=(\Delta _1, \ldots , \Delta _k)$
be the $2^1$-skeleton of $w$ and let $w_i=w|_{\Delta _i}$ and let $\eta w_i$ be their identification 
with paths in $W_1\cup V_1$ as defined in \eqnu{eta}.
Let $\sum _u$  denote the sum taken over $u\in \hat{W}_N$ satisfying $\sigma _0(u)=\frac{1}{2}\sigma _1(w)$, 
which consists of $\Delta _1, \ldots , \Delta _k$ scaled by $1/2$. 

Thus, we have  
\begin{eqnarray*}
\hat{P}_{N+1}[\ w\ ]&=&P_{N+1}[\ v\ : \ Lv=w\ ]
\\
&=&
\sum _{u} P_{N+1}[\  Lv=w,\ \frac{1}{2} \hat{Q}_1v=u\ ]
\\
&=&
\sum _{u} P_{N+1}[\  Lv=w\ |\ \frac{1}{2} \hat{Q}_1v=u\ ]
 \ P_{N+1}[\ \frac{1}{2} \hat{Q}_1v=u\ ]
\\
&=&
\sum _{u} P_{N+1}[\  Lv=w\ |\ \frac{1}{2} \hat{Q}_1v=u\ ]
\ \hat{P}_{N}[\ u\ ]
\\
&=&
\sum _{u} P_{N+1}[\ \eta (Lv|_{\Delta _i})=\eta w_i , i=1, \ldots , k\ |\ \frac{1}{2} \hat{Q}_1v=u\ ]
\ \hat{P}_{N}[\ u\ ]
\\
&=&
\sum _{u} \ (\prod _{i=1}^{k} \hat{P}^*_{1}[\ \eta w_i\ ])
\ \hat{P}_{N}[\ u\ ],
\end{eqnarray*}
where  $\hat{P}^*_{1}=\hat{P}_1$ if $\Delta _i$ is of Type 1, and 
$\hat{P}^*_{1}=\hat{P}'_1$ if $\Delta _i$ is of Type 2.

Since taking the sum over $w \in \hat{W}_{N+1}$ means taking the sum over 
all $u\in \hat{W}_N$ and finer structures in each $\Delta \in \sigma _1(w)$, we have  
\begin{eqnarray*}
\Phi _{N+1}(x,y)&=&\sum _{w\in \hat{W}_{N+1} } \hat{P}_{N+1}(w)x^{s_1(w)}y^{s_2(w)}
\\
&=&
\sum _{u\in \hat{W}_{N }} \sum _{\eta w_1 \in \hat{W}^*_1}\cdots \sum _{\eta w_k \in \hat{W}^*_1}
\ (\prod _{i=1}^{k} \hat{P}^*_{1}[\ \eta w_i\ ])
\ \hat{P}_{N}[\ u\ ]\ x^{s_1(w_1)+\cdots +s_1(w_k)}y^{s_2(w_1)+\cdots +s_2(w_k)}
\\
&=&
\sum _{u\in \hat{W}_{N }}\  \hat{P}_{N}[\ u\ ]\ \prod _{i=1}^{k} \ (\sum _{ w_i \in \hat{W}^*_1} \hat{P}^*_{1}[\  w_i\ ]
\ x^{s_1(w_i)}y^{s_2(w_i)}\ )
\\
&=&
\sum _{u\in \hat{W}_{N }}  \hat{P}_{N}[\ u\ ]\ \Phi (x,y)^{s_1(u)}\Theta (x,y) ^{s_2(u)}
\\
&=&
\Phi _{N}( \Phi (x,y), \Theta (x,y)).
\end{eqnarray*}

The calculations for 
$\hat{P}'_{N+1}$ and $\Theta _{N+1}(x,y)$ are similar.
\QED
\prfe

Define the mean matrix by
\eqnb
\eqna{matrix}
{\bf M}=\left [
\begin{array}{cc}
\frac{\partial}{\partial x}\Phi (1,1) 
& \frac{\partial}{\partial y}\Phi (1,1) 
\\
\frac{\partial}{\partial x}\Theta (1,1) 
& \frac{\partial}{\partial y}\Theta (1,1)
\end{array}
\right ]
=
\left [
\begin{array}{cc}
\frac{9}{5} & \frac{2}{15}\\
\frac{26}{15} & \frac{13}{15}
\end{array}
\right ].
\eqne
It is a strictly positive matrix, and the 
larger eigenvalue is 
\[\lambda =\frac{1}{15}(20+\sqrt{205})=2.2878 \ldots .\]

The loop-erasing procedure together with \prpu{q2} leads to
\parr
\prpb
\prpa{exittime} 
Let $M\leq N$.
Conditioned on $\sigma _M(LX_N)=(\Delta _1, \ldots , \Delta _k)$ and the types of 
each element of the skeleton, the traverse times of the triangles 
\[T_i^{ex, M}(LX_N)-T_{i-1}^{ex, M}(LX_N), \ \ i=1,2, \ldots , k\]
are independent.  Each of them  has the same distribution as either 
$T^{ex, N-M}_1(LX_{N-M})$ 
\parr
or $T^{ex, N-M}_1(LX'_{N-M})$, according to whether 
$\Delta _i$ is of Type 1 or Type 2. 
\prpe


\thmb
\thma{steps}

As $N \rightarrow \infty $, 
$\lambda ^{-N} \ell (LX_N)$ 
converges in law to an integrable random variable $W'$, 
with a positive probability density. 


\thme

We shall prove the above theorem in \secu{3}, using coupling argument. 
\thmu{steps} suggests that the displacement exponent for the loop-erased random walk on the pre-{\sg} 
 is $\log \lambda /\log 2$, 
in the sense that the average number of steps it takes to cover the distance of $2^N$
is of order $\lambda ^N$.
In other words, if we write $m=2^N$, it takes $m^{\log \lambda /\log 2}$ steps to travel a distance of $m$
from the origin.  This value is equal to that obtained by Shinoda \cite{shinoda} who defined a loop-erased walk through 
uniform spanning trees. 


\section{Scaling limit of the loop-erased random walks.}
\seca{3}

\vspace{0.5cm}\parr
\subsection{Paths on the {\sg}.}

In this section we investigate the limit of the loop-erased random walk as the lattice spacing (edge 
length) tends to $0$.
First we define the (finite) {\sg}. Since it will be easier to deal with continuous functions from the beginning, 
we regard $F_0$ as a closed subset of ${\mathbb R}^2$ made up of all the points on its 
edges.  Let $\Delta _1 $ be the  closed (filled) triangle in ${\cal T}_0$ whose vertices are $O, a_0$ and $b_0$, and    
$\Delta _2$ be its reflection with regard to the $y$-axis, and 
let $F^N=2^{-N}F_0 \cap ( \Delta _1 \cup \Delta _2)$ (Fig 5). 
We define the {\bf {\sg}}  by 
$F=cl(\cup_{N=0}^{\infty } F^{N})$, where $cl$ denotes closure.  
We define the sets of vertices by $G^N=2^{-N}G_0 \cap ( \Delta _1 \cup \Delta _2) $. 
\\
\input{Fig5.tex}\\
\\
\hspace{3cm} {\bf Fig. 5}
\vspace{0.5cm}

Let 
\[
C=\{w\in C([0,\infty)\rightarrow F) \ :\  w(0)=O,\ 
\lim_{t\rightarrow\infty} w(t)=a_0\}\,. 
\]
$C$ is a complete separable 
metric space with the metric
\[
d(u,v)=\sup_{t\in[0,\infty)}\, |u(t)-v(t)|\,,\ 
u,v\in C,
\]
where  $|x-y|$, $x, y \in \twreals$,  denotes the Euclidean distance.
Throughout this section, 
for $w \in \bigcup _{N=1}^{\infty }W_N $, we let 
\[w(t)=a_N,\ \ \ t\geq \ell (w),\]
and 
interpolate all the paths linearly, 
\[w(t)=(i+1-t)w(i)+(t-i)w(i+1),\ \   i \leq t <i+1,\ \ i=0,1,2, \cdots \]
so that we can regard $w$ as a continuous function on $[0, \infty)$. 

Let 
\[W^N=2^{-N}W_N=\{2^{-N} w \ :\ w\in W_N \},\ \ \hat{W}^N=2^{-N}\hat{W}_N,\]
where all the paths in $W^N$ are understood  to have been linearly interpolated. 
In the following we shall use this identification modulo linear interpolation.  
Thus, $W^N$ and $\hat{W}^N$ are subsets of $C$. 
For $w \in W^N$, let $\tilde{\ell}(w)=\ell(2^Nw)$.
Namely, $\tilde{\ell } (w)$ is the number of $2^{-N}$-sized `steps' the path $w$ takes to 
get to $a_0$. 

We define hitting times, coarse-graining, exit times and skeletons similarly to \secu{2}, 
but with $G_M$ replaced by $G^M$. 
Namely,  for $w\in C$ we define a sequence $\{T^M_i(w)\}_{i=0}^{m}$ of the hitting times of $G^M$, 
as follows:
$T_{0}^{M}(w)=0$, and for $i\geq 1$, let
$T_{i}^{M}(w)=\inf \{j>T_{i-1}^{M}(w) :\  w(j)\in G^{M}\setminus
\{w(T_{i-1}^{M}(w))\}\}$.  $m$ is the smallest integer such that $T_{m+1}^M(w)=\infty $.
For the hitting times we are using the same notation but we hope no 
confusion arises.
For $N \in \integers_{+}$, we define a coarse-graining  map
$Q^{N}:C \rightarrow C$ by $(Q^{N}w)(i)=w(T_{i}^{N}(w))$ 
for $i=0,1,2,\ldots, m$,  and by using linear
interpolation 
\[ (Q^{N}w)(t) = \left\{ \begin{array}{ll} (i+1-t)\ (Q^{N}w)(i)+(t-i)\
(Q^{N}w)(i+1), & i \leq t < i+1,\ i=0,1,2,\ldots, m-1, \\
  a_0, & t \geq m. \end{array} \right. \]
Notice that 
\eqnb \eqna{QQ}
Q^{M}\circ Q^{N}=Q^{M}, \ \ \ \mbox{  if  } \ \  M \leq N\eqne
holds.

Since we have defined the hitting times for every 
$w \in C$, we can define its  $2^{-M}$-{\bf skeleton}, 
$\sigma ^M(w)$  (a sequence of $2^{-M}$-triangles 
$w$ passes through) and the {\bf exit times}  $\{ T_i^{ex, M}\}$
similarly to their counterparts in \secu{2}.  
To define the loop erasing operator, recall that if $w\in W^N$, then $2^Nw \in W_N$ and $L(2^Nw )\in \hat{W}_N$ 
(modulo linear interpolation). 
Thus we define loop erasure $\tilde{L} :\ \bigcup _{N=0}^{\infty }W^N \to \bigcup _{N=0}^{\infty }\hat{W}^N $
by letting  $\tilde{L}w=2^{-N}L(2^Nw) \in \hat{W}^N$ for $w\in W^N$, $N\in {\mathbb Z}_+$, and  
we define also  
$\hat{Q}^Mw = 2^{-N}\hat{Q}_M(2^N w)\in \hat{W}^M$ for $M\leq N$. 
The only differences from the previous section are that paths are continuous (by linear interpolation) and 
confined in two neighboring unit triangles,  and that we erase loops 
from  $2^{-1}$-scale down.  For each $N\in {\mathbb Z}_+$, 
let $P^N$ be the random walk path measure on $F^N$  (a probability measure on $C$ supported on $W^N$), 
namely $P^N[w]=P_N[2^Nw]$, for $w\in W^N$.  In the following, we will focus on $P^N$. 
$V_N$'s and $P'_N$'s introduced in the previous section have played 
auxiliary roles.

\vspace{0.5cm}\parr
\subsection{The scaling limit.}

We consider random walks (linearly interporated version) on $G^N$, $N \in {\mathbb Z}_+$, starting at $O$ and stopped at 
$a_0$.

Let 
\[\Omega '=\{  \omega=(w_0, w_1, w_2, \cdots )\ :\ w_0\in \hat{W}^0,\ w_N\in \hat{W}^N,\  
w_N \triangleright w_{N+1}, \ N \in {\mathbb N} \},
\]
where $w_N  \triangleright w_{N+1}$ means that there exists a $v \in W^{N+1}$ such 
that $Q^Nv=w_N$ and $\hat{Q}^{N+1}v=w_{N+1}$. 
Define the projection onto the first $N+1$ elements by 
\[\pi _N \omega =(w_0, w_1, \ldots , w_N),\]
and a probability measure on $\pi _N\Omega '$ by 
\[\hat{P}^N[(w_0, w_1, \ldots , w_N)]=P^N[\ v :\ \hat{Q}^i v=w_i , i=0, \ldots , N\ ]\]
The following consistency condition is a direct consequence of the  loop-erasing procedure: 

\eqnb
\eqna{consistency}
\hat{P}^N [(w_0, w_1, \ldots , w_N) ]=\sum _{u} \hat{P}^{N+1} [(w_0, w_1, \ldots , w_N, u ) ],
\eqne 
where the sum is taken over all possible $u$ such that $w_N\triangleright u$.

By virtue of \eqnu{consistency} and Kolmogorov's extension theorem for a projective limit, 
there is a probability measure $\hat{P}$ on $\Omega _0=
C^{\nintegers}=C \times C \times \cdots \ $  such that
\[\hat{P}[\ \Omega '\ ]=1.\]
\[  \hat{P}\circ \pi _N^{-1} =\hat{P}^{N}, \ N\in {\mathbb Z}_+.\]

Let $Y^N :\ \Omega ' \to C$ be the projection to the $N$-th component. 
We regard  $Y^N$ as an $F$-valued process $Y^N(\omega , t )$ on 
$(\Omega _0,  {\cal B}, \hat{P})$, 
where ${\cal B}$ is the Borel
algebra on $\Omega _0$ generated by the cylinder sets.


For $w \in C$ and $j=1,2$, denote by $S_j^M(w)$ the number of $2^{-M}$-triangles of Type $j$ in $\sigma ^M(w)$, 
namely, $S^M_j(w)=\sharp \{i  :\ \Delta _i \mbox{ is of Type }j \}$,  
and let ${\bf S}^M(w)=(S_1^M(w), S_2^M(w))$.
If $w \in W^N $ for some $N$, then $\tilde{\ell} (w)=S_1^N(w)+2 S_2^N(w)$.

Let ${\bf S}=(S_1,S_2)$ and  ${\bf S}'=(S'_1,S'_2)$ be ${\mathbb Z}_+$-valued random variables on 
$(\Omega _0,  {\cal B}, \hat{P})$ 
with the same distributions as those of $(s_1, s_2)$ under $\hat{P}_1$ and  under $\hat{P}'_1$, respectively.  
$(s_1, s_2)$ has been defined in 2.4 together with the generating functions.


\prpb
\prpa{brpr}
Fix arbitrarily $v \in \hat{W}^M$, and let $\sigma ^M(v)=(\Delta _1, \ldots , \Delta _k)$.
For each $i$, $1\leq i \leq k$,   
under the conditional probability $\hat{P}[\ \ \cdot \ \ |Y^M=v ]$, 
$\{{\bf S}^{M+N}(Y^{M+N}|_{\Delta _i}),\ N=0,1,2, \cdots  \}$ is a two-type 
supercritical branching process, with the types of children corresponding  to the types of triangles.
The offspring distributions born from a Type 1 triangle 
and from a Type 2 triangle are equal to those of 
${\bf S}$ and ${\bf S}'$, respectively.  
If $\Delta _i$ is of Type 1, the process initiates in state $(1,0)$, and  
if $\Delta _i$ is of Type 2,  in state $(0,1)$.

\begin{enumerate} 
\item[(1)]
The generating functions for the offspring distributions are 
\[g_1(x, y) \defd \hat{E}[\ x^{S_1}y^{S_2}\  ]=\Phi (x,y),\]
\[g_2(x, y) \defd \hat{E} [\ x^{S'_1}y^{S'_2}\ ]=\Theta  (x,y),\]
where $ \hat{E}$ is an expectation with regard to $\hat{P}$. 
\item[(2)] 
The mean matrix ${\bf M}$ is given by \eqnu{matrix} in \secu{2}. 
It is  strictly positive and its eigenvalues are $\lambda
=\frac{1}{15}(20+\sqrt{205})=2.2878 \ldots $ and $\lambda '
=\frac{1}{15}(20-\sqrt{205})=0.3788 \ldots $.  We have
\[\hat{E}[\ {\bf S}^{M+N}(Y^{M+N}|_{\Delta _i})\  |\ Y^M=v\  ]=
{\bf S}^M(v|_{\Delta _i}) {\bf M}^N.\]
\item[(3)]
$\hat{P}[S_1+S_2 \geq 2]=\hat{P}[S_1'+S'_2 \geq 2]=1$ (non-singularity).
\item[(4)]
\[\hat{E}[\ S_i\log  S_i\ ]<\infty ,\ \ \hat{E}[\ S'_i\log  S'_i\ ]<\infty ,
\ i=1,2.\]
\end{enumerate}
\prpe

\prpu{brpr} suggests that we should consider $F$-valued processes
with time appropriately scaled.
Thus, we introduce a time-scale transformation $U_{N}(\alpha ):C \rightarrow 
C, \ \alpha \in (0,\infty),\ n \in \nintegers$.  
For $w \in C$, define 
\[(U_{N}(\alpha )w)(t) \defd w(\alpha ^{N}t),\]
and consider the processes 
\[X^N=U_{N}(\lambda )Y^N, \ \ N \in {\mathbb Z}_+ .\]


\prpb
\prpa{X}
\[\sigma ^M(X^N)=\sigma ^M(X^M)=\sigma ^M(Y^M), \ \ M \leq N ,\ \ \mbox{ a.s.}\]
\parr
In particular,
\eqnb \eqna{YY}
X^N(T_i^{ex, M}(X^N))=X^M(T_i^{ex, M}(X^M))=Y^M(T_i^{ex, M}(Y^M)),\ \ M\leq N, \ \ \mbox{ a.s.}
\eqne
\prpe
\vspace{0.5cm}\par
Note that 
if $\sigma ^M(X^N)=(\Delta _1, \cdots , \Delta _k)$, then
\[T_j^{ex, M}(X^N)=\lambda ^{-N}\sum _{i=1}^{j}(S_1^N(X^N|\Delta _i)+2S_2^N(X^N|\Delta _i)), \ \ 
1\leq j \leq k.\]

\prpu{brpr} combined with the convergence theorem for supercritical branching processes 
(see  \cite{AN}, Chapter V )
leads to the following proposition.


Let ${\bf u}=(u_1, u_2)$ and ${\bf v}=(v_1, v_2)$ be the right and left positive eigenvectors associated with $\lambda $ 
such that $|{\bf u}|= |{\bf v}|=1$.
\prpb
\prpa{supbr}
Fix arbitrarily $v \in \hat{W}^M$, and let $\sigma ^M(v)=(\Delta _1, \ldots , \Delta _k)$.
For each $i$, $1\leq i \leq k$,   
under the conditional probability $\hat{P}[\ \ \cdot \ \ |Y^M=v ]$, 
we have the following.
\itmb
\item[(1)]
For each $i \in \{1, \cdots , k\}$, 
$\{\lambda ^{-(M+N)}{\bf S}^{M+N}(X^{M+N}|_{\Delta _i}),\ N=0,1,2, \ldots \}$ converges a.s. 
as $N \rightarrow \infty$
to a ${\mathbb R}^2$-valued random variable ${\bf S}^{*M,i}=(S_1^{*M,i}, S_2^{*M,i})$.
\item[(2)]
$\{ {\bf S}^{*M,i}$, $i=1, \cdots , k\}$ are independent.
\item[(3)]
There are random variables $B_1$ and $B_2$ such that 
${\bf S}^{*M,i}$ is equal in distribution to $\lambda ^{-M}B_1{\bf v}$ 
if $\Delta _i$ is of Type 1, and equal in distribution to 
$\lambda ^{-M}B_2{\bf v}$ 
if $\Delta _i$ is of Type 2. 

\item[(4)] 
\[\hat{P}[B_i>0]=1,\ \    
\hat{E}[B_i]=u_i,\ \ i=1,2.\] 
$B_1$ and $B_2$ have strictly positive probability density functions.

\item[(5)]
 
The Laplace transform of  $B_i$, $i=1,2$    
\[\phi _i(t)=\hat{E}[\exp ( tB_i)]\]
are entire functions on ${\mathbb C}$ and are the unique solution to 
\[\phi _1 (\lambda  t)=\Phi( \phi _1( t), \phi _2 ( t)),\  \  
\phi _2 (\lambda  t)=\Theta ( \phi _1( t), \phi _2( t)),\
\phi _1(0)=\phi _2(0)=1.
\]

\itme
\prpe

To be precise, 
(1)--(4) in \prpu{supbr} are the straightforward consequences of general limit theorems
for superbranching processes (Theorem 1 and Theorem 2 in V.6 of \cite{AN}).
$\hat{P}[B_i>0]=1$ is a consequence of $\Phi$ and $\Theta$ having no terms with 
degree smaller than $2$.  
For the exsistence of the Laplace transform on the entire ${\mathbb C}$, 
we need careful study of the recursions.  We omit the details here, since they are 
similar to those in \cite{hhk2}.

Let $T^{*M}_{i}=\sum_{j=1}^{i}(S_1^{*M,j}+2S_2^{*M,j})$.  Then $\dsp \lim _{N \to \infty}
T_j^{ex, M}(X^N)=T_j^{*M}$.
By virtue of \prpu{X} and \prpu{supbr}, we can prove the almost sure uniform convergence
for $X^N$.

\thmb
\thma{asconv}
$X^N$ converges uniformly in $t$ a.s. as $N \rightarrow \infty$ to a 
continuous process $X$. 
\thme

\prfb
Choose $\omega \in \Omega'$ such that for all 
$M \in {\mathbb Z}_+$ the following holds: 
$Y^M \in \hat{W}^M$, 
$\dsp \lim _{N \to \infty} T_i^{ex, M}(X^N)=T_i^{*M}$ exists 
and $T_i^{*M}-T_{i-1}^{*M}>0$
for all  $1\leq i \leq k_M$, where $k_M$ denotes the number of triangles in $\sigma ^M(Y^M)$.  
Let $R=T_1^{*0}+ \varepsilon$,
where $\varepsilon >0$ is arbitrary.  
It suffices to show that $X^N(\omega ,t)$ converges uniformly 
in $t\in [0,R]$.  In fact, if $t>R$, $X^N(t)=a_0$ for a large enough $N$.

Fix $M \geq 0$.   Let $k=k_M$. 
By expressing the arrival time at $a_0$ as the sum of traversing times of $2^{-M}$-triangles, 
we have 
$T^{ex, M}_{ k}(X^N)=T_1^{ex, 0}(X^N)$ a.s.  
Letting $N \rightarrow \infty$, 
we have $T^{*M}_{k}=T_{1}^{*0}$ a.s.   

The choice of $\omega$ shows that 
there exists an $N_1=N_1(\omega ) \in \nintegers$
such that
\eqnb
\eqna{max}
\max_{1 \leq i \leq k}|T_i^{ex, M} (X^N)-T_i^{*M}| \leq 
\min_{1 \leq i \leq k}
(T_i^{*M}-T_{i-1}^{*M})
,\eqne   
and \[|T_{k}^{ex, M}(X^N)-T_{k}^{*M}|<\varepsilon ,\]
for $N \geq N_1$.

If $0\le t< T^{*M}_{k}$, then choose
$j \in \{1, \cdots , k\}$ such that 
$T^{*M}_{j-1} \leq t < T^{*M}_{j}$.

Then \eqnu{max} implies  that $T_{j-2}^{ex, M}(X^N)\leq t \leq T_{j+1}^{ex, M}(X^N)$,
for $N \geq N_1$.  Since \prpu{X} shows 
\eqnb \eqna{XX}
X^N(T_j^{ex, M}(X^N))=X^M(T_j^{ex, M}(X^M)), \eqne
for all $N$ with $N\geq M$, we have 
\[|X^N(T_j^{ex, M}(X^N))-X^N(t)| \leq 3 \cdot 2^{-M}. \]
Otherwise, if $ T^{*M}_{k}\le t \le  T^{*M}_{k}+\varepsilon =R$, then 
let $j=k$.
Since $T_{k-1}^{ex, M}(X^N) \leq t$,
\[|X^N(T_j^{ex, M}(X^N))-X^N(t)| \leq 2 \cdot 2^{-M}. \]

Therefore, if  $N,N' \geq N_1$, then for any $t \in [0,R]$,
\begin{eqnarray*}
\lefteqn
{|X^N(t)-X^{N'}(t)|
}\\
&\leq&
|X^N(T_j^{ex, M}(X^N))-X^{N}(t)|+|X^{N'}(T_j^{ex, M}(X^{N'}))-X^{N'}(t)|\\
&&
+|X^N(T_j^{ex, M}(X^N))-X^{N'}(T_j^{ex, M}(X^{N'}))|
\\
& \leq &
6 \cdot 2^{-M},
\end{eqnarray*}
where the third term in the middle part is shown to be $0$ by \eqnu{XX}.
Since $M$ is arbitrary, we have the uniform convergence.

\QED
\prfe

\thmb
\thma{SA}
$X$ is almost surely self-avoiding. 
The Hausdorff dimension of the path $X([0, \infty ))$ is almost surely equal to $\log \lambda /\log 2$.
\thme

The uniform convergence of $X^N$, which is self-avoiding,  to $X$ implies that 
the probability of the event that 
there exist $t_1$, $t_2$ and $t_3$ with $t_1<t_2<t_3$ such that 
$X(t_1)=X(t_3)$,\  \ $ X(t_2) \neq X(t_1)$ 
is zero, and the existence of the Laplace transforms $\hat{E}[\exp (t_0B_i)]$, $i=1,2$ for some $t_0>0$ 
guarantees that the probability that there exist $t_1, t_2>0$ such that 
$X(t)=X(t_1)$ for all $t$, $t_1\leq t\leq t_1+t_2$ is zero. 
We omit the proof here since they are similar to that in  \cite{hhk2}.
To calculate the Hausdorff dimension, we regard the path as a multi-type random 
fractal. The proof is similar to that in  \cite{KH}.  

\vspace{0.5cm}\par
Since 
$\lambda ^{-N}\ell (LX_N)$ 
in \thmu{steps} has the same distribution as 
$\lambda ^{-N}(S_1^N(X^N)+2S_2^N(X^N))$, 
\thmu{steps} follows immediately from \prpu{supbr}, with $W'$ equal in distribution to $(v_1+2v_2)B_1$.

\section{Conclusion}
\seca{4}

We proposed a model of loop-erased random walks on the finite pre-{\sg}.  
Our loop-erasing procedure is based on a `larger-scale-loops-first' rule, which enables us to 
obtain exact recursion relations.  
First, we proved the existence of the scaling limit.  Then, we made use of the 
tools that have been developped for the study of self-avoiding walks on the pre-{\sg} 
to prove that the path of the limiting process is almost surely self-avoiding, 
while having Hausdorff dimension strictly greater than $1$.  
Our path Hausdorff dimension is consistent with the results in \cite{DD} and 
\cite{shinoda}, thus we conjecture that our model is in the same universality 
class as theirs.


\begin{thebibliography}{2007}

\bibitem{AN} K.~B.~Athreya, P.~E.~Ney, 
{\it Branching processes,}
Springer, 1972.

\bibitem{BP}
M.T.~Barlow, E.A.~Perkins, 
\textit{Brownian motion on the Sierpinski gasket}, 
Probab. Theory Relat. Fields \ \textbf{79} (1988) 543--623

\bibitem{DD}
D.~Dhar, A.~Dhar
\textit{Distribution of sizes of erased loops for loop-erased random walks}
Physical Review E,  \ \textbf{55} (1997) R2093--2096

\bibitem{falconer}
K.~Falconer, \textit{Fractal geometry}, 2nd ed. Wiley, 2003

\bibitem{HHH} B.~Hambly, K.~Hattori, T.~Hattori,
{\it Self-repelling walk on the {\sg},}
Probab.\ Theory Relat.\ Fields \ {\bf 124} (2002) 1--25.


\bibitem{hhk1}
K.~Hattori, T.~Hattori, S.~Kusuoka,
\textit{Self-avoiding paths on the pre-Sierpinski gasket,} 
Probab. Theory Relat. Fields \ \textbf{84} (1990) 1--26

\bibitem{hh}
K.~Hattori, T.~Hattori, 
\textit{Self-avoiding process on the Sierpinski gasket,} 
Probab. Theory Relat. Fields \ \textbf{88} (1991) 405--428

\bibitem{hk}
 T.~Hattori, S.~Kusuoka,
\textit{The exponent for mean square displacement of self-avoiding random walk  on Sierpinski gasket,} 
Probab. Theory Relat. Fields \ \textbf{93} (1992) 273--284.

\bibitem{hhk2}
K.~Hattori, T.~Hattori, S.~Kusuoka,
\textit{Self-avoiding paths on the three-dimensional Sierpinski gasket,}
Publication of RIMS \ \textbf{29} (1993) 455--509

\bibitem{KH}
K.~Hattori, 
\textit{Fractal geometry of self-avoiding processes}, J. Math. Sci. Univ.
 Tokyo,  \ \textbf{3} (1996) 379--397


\bibitem{Kozma}
G.~Kozma, 
\textit{The scaling limit of loop-erased random walk in three dimensions}, Acta Math. \ \textbf{1} (2007) 29--152


\bibitem{lawler}
G.F.~Lawler, 
\textit{Intersection of random walks},  \ Birkh\"auser, 1991

\bibitem{LSW}
G.F.~Lawler, O.~Schramm, W.~Werner,
\textit{Conformal invariance of planar loop-erased random walks and uniform spanning 
trees}, Ann. Probab. 
\ \textbf{32} (2004) 939--995.

\bibitem{ms}
N.~Madras, G.~Slade,
\textit{The self-avoiding walk} \ Birkh\"auser, 1993

\bibitem{schramm}
O.~Schramm
\textit{Scaling limits of loop-erased random walks and uniform spanning trees},
Israel Journal of Mathematics, \ \textbf{118} (2000), 221--288.

\bibitem{shinoda}
M.~Shinoda, 
\textit{Uniform spanning tree measure on the pre-\sg}, (in Japanese)
unpublished


\end{thebibliography}
\end{document}